\documentclass[MSNbibl,number,citesort,seceqn,dvips]{arxbj}

\aid{0}
\volume{20}
\issue{3}
\pubyear{2014}
\firstpage{1432}
\lastpage{1453}
\doi{10.3150/13-BEJ528} 

\makeatletter
\newcommand{\mrmd}{\,\mathrm{d}}
\newcommand{\mrmdd}{\mathrm{d}}

\newcommand{\rrvert}{\vert}
\newcommand{\llvert}{\vert}

\newcommand{\cal}{\mathcal}

\newtheorem{theor}{Theorem}[section]
\newproclaim{assumption}{Assumption}

\newcommand{\R}{\mathbb R}
\newcommand{\N}{\mathbb N}

\newcommand{\pr}{^{\prime}}

\newcommand{\thetab}{{\bolds\theta}}
\newcommand{\varthetab}{{\bolds\vartheta}}

\makeatother

\begin{document}
\begin{frontmatter}

\title{Skew-symmetric distributions and Fisher information: The double
sin of\\  the skew-normal}
\runtitle{The double sin of the skew-normal}

\begin{aug}
\author[1]{\inits{M.}\fnms{Marc} \snm{Hallin}\corref{}\thanksref{1}\ead[label=e1]{mhallin@ulb.ac.be}} \and
\author[2]{\inits{C.}\fnms{Christophe} \snm{Ley}\thanksref{2}\ead[label=e2]{chrisley@ulb.ac.be}}
\runauthor{M. Hallin and C. Ley} 
\address[1]{ECARES, Universit\'{e} Libre de Bruxelles,
CP114, 50 Ave. F.D. Roosevelt,
1050 Brussels, Belgium, and ORFE, Princeton University,
Sherrerd Hall, Princeton, NJ 08544, USA.\\
\printead{e1}}
\address[2]{ECARES and D\'epartement de Math\'ematique, Universit\'{e}
Libre de Bruxelles, CP210,
Boulevard du Triomphe, 1050 Brussels, Belgium. \printead{e2}}
\end{aug}

\received{\smonth{9} \syear{2012}}
\revised{\smonth{4} \syear{2013}}

%
\begin{abstract}
Hallin and Ley [\textit{Bernoulli} \textbf{18} (2012) 747--763]
investigate and fully characterize the Fisher singularity phenomenon in
univariate and multivariate families of \textit{skew-symmetric
distributions}. This paper proposes a refined analysis of the
(univariate) problem, showing that singularity can be more or less
severe, inducing $n^{1/4}$ (``simple singularity''), $n^{1/6}$
(``double singularity''), or $n^{1/8}$ (``triple singularity'')
consistency rates for the skewness parameter. We show, however, that
simple singularity (yielding $n^{1/4}$ consistency rates), if any
singularity at all, is the rule, in the sense that double and triple
singularities are possible for \textit{generalized skew-normal} families
only. We also show that higher-order singularities, leading to
worse-than-$n^{1/8}$ rates, cannot occur. Depending on the degree of
the singularity, our analysis also suggests a simple reparametrization
that offers an alternative to the so-called \textit{centred
parametrization} proposed, in the particular case of skew-normal and
skew-$t$ families, by Azzalini [\textit{Scand. J. Stat.} \textbf{12}
(1985) 171--178], Arellano-Valle and Azzalini [\textit{J. Multivariate
Anal.} \textbf{113} (2013) 73--90], and DiCiccio and Monti
[\textit{Quaderni di Statistica} \textbf{13} (2011) 1--21],
respectively.
\end{abstract}

%

\end{frontmatter}

\section{Introduction}

The \emph{skew-symmetric} families, originally proposed in
Azzalini and Capitanio \cite{AC03} and Wang, Boyer and Genton \cite
{WBG04}, are, in their
univariate version, parametric families of probability density
functions (p.d.f.s) of the form
%
%
\begin{equation}
\label{HLSS} x\mapsto f_{\varthetab}^\Pi(x):=2 \sigma^{-1}f
\bigl(\sigma^{-1}(x-\mu)\bigr)\Pi\bigl(\sigma^{-1}(x-\mu),
\delta\bigr),\qquad x\in\R,
\end{equation}
\begin{enumerate}[(c)]
\item[(a)] $\varthetab=(\mu,\sigma,\delta)\pr$, with $\mu\in\R
$ a \textit{location parameter} and $\sigma\in\R_0^+$ a \textit
{scale parameter}, while $\delta\in\R$ plays the role of a \textit
{skewness parameter};
\item[(b)] $f\dvtx \R\rightarrow\R^+_0$, the \emph{symmetric kernel},
is a nonvanishing symmetric p.d.f. (such that, for any $ z\in\R$, $0\neq
f(-z)=f(z)$), and
\item[(c)] $\Pi\dvtx\R\times\R\rightarrow[0,1]$ is a \emph{skewing
function}, that is, satisfies
%
%
\begin{equation}
\label{c1}\Pi(-z,\delta)+\Pi(z,\delta)=1,\qquad z,\delta\in\R\quad
\mbox{and}\quad
\Pi(z,0)=1/2, \qquad z\in\R,
\end{equation}
and, in case $(z,\delta)\mapsto\Pi(z,\delta)$ admits a derivative
of order $s$ at $\delta=0$ for all $z\in\R$,
%
%
\begin{eqnarray}
\label{c2}\partial^s_z\Pi(z,\delta)\vert_{\delta
=0}&=&0,\qquad
z\in\R\quad\mbox{and}\nonumber\\[-8pt]\\[-8pt]
\mbox{for $s$ even}\qquad \partial^s_\delta\Pi(z,
\delta)\vert_{\delta=0}&=&0,\qquad z\in\R.\nonumber
\end{eqnarray}
\end{enumerate}

The assumption of a nonvanishing kernel $f$ is not essential, and has
been added, as in Hallin and Ley \cite{HL12}, in order to ease calculations
and avoid trivial complications. Condition (\ref{c1}) is classical;
the less classical condition (\ref{c2}) is justified by the analogy
with skewing functions of the form $\Pi(z,\delta)=\Pi(\delta z)$, by
far the most common ones: if $\Pi$ is $s$ times continuously
differentiable, $\partial^s_z\Pi(\delta z)=\delta^s(\partial^s\Pi
)(\delta z)$ obviously vanishes at $\delta=0$. Similarly, the fact
that $\Pi(-y)+\Pi(y)=1$ implies that $\partial^s\Pi(\delta z)$
cancels at $\delta=0$ for even values of $s$. All skewing functions
considered in the literature, as well as those appearing in the
examples developed in this paper and in Hallin and Ley \cite{HL12}, satisfy
(\ref{c2}). Further comments on  skewing functions of the form $\Pi
(z,\delta)=\Pi(\delta z)$ can be found in Section~\ref{fc}.

The \emph{skew-normal} family of Azzalini \cite{A85}, for which the
symmetric kernel $f$ is the standard Gaussian p.d.f. $\phi$ and the
skewing function $\Pi(z,\delta)=\Phi(\delta z)$ for $\Phi$ the
standard Gaussian cumulative distribution function (c.d.f.), is by far the
oldest and most popular example of such a skew-symmetric family;
varying $f$ and $\Pi$, however, yields a virtually infinite number of
them. Traditional examples include
the \textit{skew-exponential power distributions} of Azzalini \cite{A86},
the \textit{skew-Cauchy distributions} of Arnold and Beaver \cite{AB00},
the \textit{skew-$t$ densities} of Azzalini and Capitanio~\cite{AC03}, or
the \textit{generalized skew-normal distributions} of Loperfido \cite{L04};
the latter result (along with arbitrary $\Pi$) from letting $f=\phi$
in (\ref{HLSS}), and play an important role in this paper.
We refer to
Genton~\cite{G04},
Azzalini \cite{A05} or
Ley \cite{L12} for background reading, details and examples.

Since the pioneering paper by Azzalini \cite{A85}, it is well known that the
scalar skew-normal distribution suffers from a Fisher information
singularity problem at $\delta=0$. More precisely, the Fisher
information matrix for the three-parameter density (\ref{HLSS}) in the
scalar skew-normal case is singular -- typically, with rank 2 instead of
3 -- in the vicinity of symmetry, that is, at $\delta=0$. Such a
singularity violates the standard assumptions for root-$n$ asymptotic
inference, and skew-normal distributions therefore may be difficult
from an inferential point of view; in particular, traditional tests of
the null hypothesis of symmetry, at first sight, seem problematic.

That singularity problem has been discussed at length in a number of
papers, among which Azzalini and Capitanio \cite{AC99}, Pewsey \cite{P00},
DiCiccio and Monti \cite{DM04}, Chiogna \cite{C05}, Azzalini and Genton
\cite{AG08}
or Ley and Paindaveine \cite{LP10a}; see Hallin and Ley \cite{HL12} for
a detailed
account. While all authors were pointing at some special status for
normal kernels, hence skew-normal distributions, Hallin and Ley \cite{HL12}
have shown that this information singularity has no special relation to
the skew-normal case, but actually originates in an unfortunate
mismatch between $f$ and $\Pi$ -- more precisely, between two
densities, the kernel $f$ and an exponential density $g_\Pi$
associated with the skewing function $\Pi$ (see Section~\ref{mismatch}).

Singularity of Fisher information results in slower consistency rates
in the estimation of the skewness parameter (at $\delta
=0$) -- equivalently, it yields slower local alternative rates (\textit
{contiguity} rates) in tests of the null hypothesis of
symmetry ($\delta=0$). That impact of singular Fisher information on
consistency/contiguity rates has been studied, in a general context,
for the particular case of a singularity of order one, by
Rotnitzky \textit{et al.} \cite{Ro00}, who unify and reinforce earlier
proposals by,
for example, Cox and Hinkley (\cite{CH74}, pages 117 and 118)
or Lee and Chesher \cite{LC86}.

The typical rate, corresponding to a ``simple singularity'', would be
$n^{1/4}$. However, it is well known (see, e.g., Chiogna \cite{C05})
that, for skew-normal distributions, that $n^{1/4}$ rate (for the
estimation of $\delta$ at $\delta= 0$) drops down to $n^{1/6}$. In
order to understand and explain this intriguing phenomenon, we pursue
and refine, in the present paper, the analysis of Fisher singularity
initiated in Hallin and Ley \cite{HL12}. We show that this deterioration
from $n^{1/4}$ to $n^{1/6}$ is explained by a ``double singularity''
property (a terminology that will become clear in the course of this
paper) -- the \textit{double sin} of the skew-normal. That $n^{1/6}$
rate in turn possibly can drop further down to~$n^{1/8}$, a case of
``triple singularity''. This, however, as we show in Theorem \ref
{result2}, is the worst case: ``fourfold singularities'' --
\textit{quadruple sins} -- yielding $n^{1/10}$ rates or worse, are excluded.

Our aim is to characterize, in the spirit of Hallin and Ley
\cite{HL12}, among all families of\linebreak[4]  univariate skew-symmetric
distributions suffering from Fisher singularity, those exhibiting that
simple/double/triple singularity phenomenon, and to show that there
exist no higher-order ones. It turns out that only Gaussian kernels can
exhibit double (a fortiori, triple) singularity. The skew-normal family
is but one example; other cases are found in the aforementioned class
of generalized skew-normal distributions (Loperfido \cite{L04}). We also
provide (in the spirit of  Rotnitzky \textit{et al.}~\cite{Ro00}) the
reparametrizations and the scores taking care of simple, double and
triple singularities, and achieving the $n^{1/4}, n^{1/6}$ and
$n^{1/8}$ consistency/contiguity rates, respectively, for $\delta$ at
$\delta=0$.

With the same objective of overcoming the Fisher singularity problem,
Azzalini \cite{A85}, Arellano-Valle and Azzalini \cite{AA08},
Arellano-Valle and Azzalini \cite{AA11} and DiCiccio and Monti
\cite{DM11}, in the univariate skew-normal, multivariate skew-normal,
and univariate skew-$t$ cases, respectively, also propose a nonlinear
reparametrization, called the \textit{centred parametrization}, which
we discuss in some detail in Section~\ref{cp}. Based on Gram--Schmidt
orthogonalization (in the space of \textit{scores}) ideas, our
reparametrizations are analytically simple, and lead to closed-form
expressions of the scores at $\delta=0$.
As far as the asymmetry parameter is concerned, and depending on the
degree of Fisher singularity, they involve
$\delta^{(1)}:=\operatorname{sign}(\delta)\delta^2$,
$\delta^{(2)}:=\delta^3$, or $\delta
^{(3)}:=\operatorname{sign}(\delta)\delta^4$, thus preserving the
interpretation of skewness as a tuning parameter of the skewing
mechanism and the characterization of the null hypothesis of symmetry
as $\delta= 0$. Finally, our reparametrizations apply to arbitrary
skew-symmetric families with singular Fisher information at $\delta=0$.

The paper is organized as follows. Section~\ref{step1} deals with the
simple singularity case, Section~\ref{step2} with double singularity.
Section~\ref{step3} analyzes the triple singularity case and shows
that higher-order ones are excluded. Examples for each type of
singularity are provided in Section~\ref{exos}. We conclude the paper
with a discussion of Azzalini's centred parametrization and its
relation to ours (Section~\ref{cp}), and some warnings (Section~\ref
{fc}) about the potential dangers of the standard skewing functions $\Pi
(z,\delta)=\Pi(\delta z)$.

\section{Simple singularity}\label{step1}

In this section, we first briefly revisit the main result of Hallin and
Ley \cite{HL12} in order, mainly, to settle the notation. We then show how
to resolve the singularity problem via an adequate reparametrization
leading, in general, to $n^{1/4}$ consistency rates for $\delta$ in
the vicinity of symmetry -- equivalently, to a reparametrization of
skewness, of the form $\delta^{(1)}:=\operatorname{sign}(\delta)\delta^2$,
recovering $n^{1/2}$ rates all over the parameter space.

\subsection{Simple singularity: A mismatch between $f$ and \texorpdfstring{$\Pi$}{Pi}}\label{mismatch}
Throughout, we consider the skew-symmetric distributions with p.d.f.
(\ref{HLSS}), along with regularity assumptions on $f$ and $\Pi$ that
will be tightened from section to section. The minimal regularity
assumptions we need are those of Hallin and Ley \cite{HL12}.

\renewcommand{\theassumption}{(A1)}
\begin{assumption}\label{assumA1}
\textup{(i)} The symmetric kernel $f$ is a \emph{standardized}
symmetric p.d.f. \textup{(ii)} The mapping $z\mapsto f(z)$ is
continuously differentiable, with derivative $\dot f$, at all
$z\in\mathbb{R}$. \textup{(iii)} Letting $\varphi_f:=-\dot{f}/f$, the
information quantities $\sigma^{-2} {\cal I}_f$ for location and
$\sigma^{-2}{\cal J}_f$ for scale, with
\[
{\cal I}_f:=\int_{-\infty}^{\infty}
\varphi^2_f(z) f(z)\mrmd  z \quad\mbox{and}\quad {\cal
J}_f:= \int_{-\infty}^{\infty}\bigl(z
\varphi_f(z)-1\bigr)^2 f(z)\mrmd  z,
\]
are finite.
\end{assumption}

\renewcommand{\theassumption}{(A2)}
\begin{assumption}\label{assumA2}
\textup{(i)} The mapping $(z,\delta)\mapsto \Pi(z,\delta)$ is
continuously differentiable at $\delta=0$ for all $z\in\mathbb{R}$;
\textup{(ii)} the derivative $ {{\partial}_{ \delta}\Pi
(z,\delta)|_{\delta=0}=:\psi(z)}$ admits a primitive $\Psi$;
\textup{(iii)} the quantity $\int_{-\infty}^\infty\psi^2(z)f(z) \mrmd
z$ is finite.
\end{assumption}

Regarding Assumption \ref{assumA1}(i), the term ``standardized'' means
that the
scale parameter (not necessarily a standard error, so that finite
second-order moments are not required) of the symmetric kernel equals
one -- an identification constraint for $\sigma$ that does not imply
any loss of generality; see Hallin and Ley \cite{HL12} for a
discussion. %
All other assumptions ensure the existence and finiteness of Fisher
information for the original parametrization.

Under Assumptions \ref{assumA1} and \ref{assumA2},
the \emph{score vector} ${\bolds\ell}_{f;\varthetab}$, at $ (\mu,
\sigma, 0)^\prime=:\varthetab_0$, takes the form
\begin{eqnarray*}
{\bolds\ell}_{f;\varthetab_0}(x)&:=& \mathrm{grad}_\varthetab\log
f_{\varthetab}^\Pi(x) \vert_{\varthetab_0}=: \bigl(
\ell_{f;\varthetab_0}^{1}(x), \ell_{f;\varthetab_0}^{2}(x),
\ell_{f;\varthetab_0}^{3}(x) \bigr)\pr
\\
&=&\pmatrix{ \sigma^{-1}\varphi_f\bigl(
\sigma^{-1}(x-\mu)\bigr)
\cr
\sigma^{-1}\bigl(
\sigma^{-1}(x-\mu)\varphi_f\bigl(\sigma^{-1}(x-
\mu)\bigr)-1\bigr)
\cr
2\psi\bigl(\sigma^{-1}(x-\mu)\bigr)},
\end{eqnarray*}
where the factor 2 in $\ell_{f;\varthetab_0}^{3}$ follows from the
fact that $\Pi(z,0)=1/2$ for all $z\in\R$. Note that the skewing
function $\Pi$ plays no role in the score functions for $\mu$ and
$\sigma$ at $\delta=0$. The resulting $3\times3$ Fisher information
matrix then exists, is finite, and takes the form
\[
{\bolds\Gamma}_{ f;\varthetab_0}:=\sigma^{-1} \int
_{-\infty
}^{\infty} {\bolds\ell}_{f;\varthetab_0}(x){\bolds
\ell}^{\prime
}_{f;\varthetab_0}(x) f\bigl(\sigma^{-1} (x-\mu)
\bigr)\mrmd x =: \pmatrix{ \gamma_{f;\varthetab_0}^{11}&0&\gamma
_{f;\varthetab_0}^{13}
\vspace*{1pt}\cr
0&\gamma_{f;\varthetab_0}^{22}&0
\vspace*{1pt}\cr
\gamma_{f;\varthetab_0}^{13}&0&
\gamma_{f;\varthetab_0}^{33}}
\]
with
\[
\gamma_{f;\varthetab_0}^{11}=\sigma^{-2} {\cal
I}_f,\qquad 
\gamma_{f;\varthetab_0}^{22}=
\sigma^{-2}{\cal J}_f,\qquad \gamma_{f;\varthetab_0}^{33}=4
\int_{-\infty}^\infty\psi^2(z)f(z)\mrmd  z
\]
and
\[
\gamma_{f;\varthetab_0}^{13}=2\sigma^{-1} \int
_{-\infty}^\infty\varphi_f(z)\psi(z)f(z)\mrmd  z.
\]
The zeroes in ${\bolds\Gamma}_{ f;\varthetab_0}$ are easily
obtained by noting that $\ell^{1}_{f;\varthetab_0}$ and $\ell
^3_{f;\varthetab_0}$ are odd functions of \mbox{$(x-\mu)$}, whereas $\ell
^2_{f;\varthetab_0}$ is even with respect to the same quantity.
Consequently, Fisher singularity only can be caused by the collinearity
of $\ell_{f;\varthetab_0}^{1}$ and $\ell_{f;\varthetab_0}^{3}$.
Starting\vspace*{1pt} from that elementary observation, Hallin and Ley \cite{HL12} show
that the family of densities (\ref{HLSS}) characterized by a couple
$(f,\Pi)$ suffers from Fisher singularity at $\delta= 0$ if and only
if the symmetric kernel $f$ belongs to the \textit{exponential family}
%
%
\begin{equation}
\label{expfam} {\cal E}_\Psi:= \biggl\{ g_a:=\exp(-a
\Psi)\Big/\int_{-\infty}^\infty\exp\bigl(-a\Psi(z)\bigr)\mrmd z {
\Big\vert} a \in{\cal A} \biggr\}
\end{equation}
with \textit{minimal sufficient statistic} $\Psi$, \textit{natural
parameter} $-a$, and \textit{natural parameter space}
\[
{\cal A}:= \biggl\{ a\in\R\mbox{ such that } \int_{-\infty}^\infty
\exp\bigl(-a\Psi(z)\bigr)\mrmd z<\infty\biggr\},
\]
yielding
%
%
\begin{equation}\label{gammas}
\gamma_{f;\varthetab_0}^{11}=\sigma^{-2}a^2 \int
_{-\infty}^\infty\psi^2(z)f(z)\mrmd  z \quad\mbox{and}\quad
\gamma_{f;\varthetab
_0}^{13}=2\sigma^{-1} a \int
_{-\infty}^\infty\psi^2(z)f(z)\mrmd  z.
\end{equation}
We refer the reader to the end of Section~2.1 in Hallin and Ley \cite{HL12}
for comments and a discussion on the existence of couples $(f,\Pi)$
such that $f\in{\cal E}_\Psi$ for given $f$ and for given $\Pi$,
respectively.

\subsection{Towards a singularity-free reparametrization:
Orthogonalization}\label{GS1}

A natural way to handle this singularity problem consists in
reparametrizing (\ref{HLSS}) in the spirit of
Rotnitzky \textit{et al.} \cite{Ro00}. Assume that $f$ and $\Pi$ are
such that $f\in{\cal
E}_\Psi$. The collinearity at $\varthetab_0$ between the score for
location and the score for skewness can be taken care of by a
Gram--Schmidt orthogonalization process applied to the three components
of ${\bolds\ell}_{f;\varthetab_0}$. This process projects, in the
$L_2$ geometry of the information matrix, the score for skewness $\ell
_{f;\varthetab_0}^{3}$ onto the subspace orthogonal (at $\varthetab
_0$) to the scores for location and scale $\ell_{f;\varthetab_0}^{1}$
and $\ell_{f;\varthetab_0}^{2}$, so that the score for skewness
becomes orthogonal to the score for location (since it is already
orthogonal to $\ell_{f;\varthetab_0}^{2}$). The resulting score for
skewness is
\[
\ell_{f;\varthetab_0 }^{3(1)}= \ell_{f;\varthetab_0 }^{3}-
\ell_{f;\varthetab_0 }^{1} {\operatorname{Cov}\bigl(\ell_{f;\varthetab
_0 }^{1},
\ell_{f;\varthetab_0 }^{3}\bigr)}/{\operatorname{Var}\bigl(\ell
_{f;\varthetab_0 }^{1}
\bigr)},
\]
while the other two scores remain unchanged: $\ell_{f;\varthetab_0
}^{1(1)} =\ell_{f;\varthetab_0 }^{1}$, $\ell_{f;\varthetab_0
}^{2(1)} =\ell_{f;\varthetab_0 }^{2}$.
As expected, in view of (\ref{gammas}), 
%
\[
\ell_{f;\varthetab_0 }^{3(1)}(x)= 2\psi\bigl(\sigma^{-1}(x-\mu)
\bigr) - \sigma^{-1} a\psi\bigl(\sigma^{-1}(x-\mu)\bigr)
\frac{2\sigma^{-1}a\int_{-\infty}^\infty\psi^2(z)f(z)\mrmd  z}{\sigma
^{-2}a^2\int_{-\infty
}^\infty\psi^2(z)f(z)\mrmd  z} = 0.
\]

This (orthogonal at $\varthetab_0$) system of scores is associated
with the reparametrization $(\mu^{(1)},\sigma^{(1)}, \delta)\pr$, where
$\mu^{(1)}=
\mu+2\delta\sigma/{a}$ and $
\sigma^{(1)}=\sigma$,
hence
%
%
\begin{eqnarray}
\label{frepar1}
&&
f_{\mu^{(1)},\sigma^{(1)}, \delta}^\Pi(x)\nonumber\\[-8pt]\\[-8pt]
&&\qquad:=
2
\bigl(\sigma^{(1)}\bigr)^{-1}f\bigl(\bigl(x-\mu^{(1)}
+2\delta\sigma^{(1)}/{a}\bigr)/\sigma^{(1)}\bigr)\Pi\bigl(
\bigl(x-\mu^{(1)} +2\delta\sigma^{(1)}/{a}\bigr)/
\sigma^{(1)},\delta\bigr).\qquad\nonumber
\end{eqnarray}
Note that this reparametrization, which
only affects the location parameter, is adopted for the family as a
whole, although the orthogonalization argument it is based on only
holds at $ (\mu^{(1)}, \sigma^{(1)}, 0)^\prime=(\mu,\sigma,0)\pr=
\varthetab_0$.

Since this reparametrization cancels, at $\delta=0$, 
the score for skewness, hence the linear term in the Taylor expansion
of the log-likelihood, second derivatives with respect to $\delta$
naturally come into the picture. To be precise, since the linear
term $\tau_3 \partial_{\delta} \log f_{\mu^{(1)},\sigma^{(1)},
\delta}^\Pi(x)\vert_{\varthetab_0}$ in the Taylor expansion of
$\log f_{{\bolds\vartheta}_0+(0,0,\tau_3)\pr}^\Pi(x)$ about
$\log f_{{\bolds\vartheta}_0}^\Pi(x)$ happens\vspace*{1pt} to be zero, the first
local approximation is provided by the quadratic term $({\tau
_3^2}/{2})\partial_\delta^2 \log f_{\mu^{(1)},\sigma^{(1)}, \delta
}^\Pi(x)\vert_{\varthetab_0} $ -- provided that second derivatives
exist. As a result, if the impact, on the log-likelihood of an i.i.d.
sample of size $n$, of a perturbation $\tau_3$ of $\delta=0$ is to
exhibit the central-limit magnitude of $n^{-1/2}$, $\tau_3$ itself has
to be of magnitude $n^{-1/4}$ only; moreover, information about its
sign is lost (a phenomenon which is also stressed by
Rotnitzky \textit{et al.} \cite{Ro00}). This is the structural reason
for slower-than-root-$n$
consistency rates (at~$\varthetab_0$) for the skewness parameter
$\delta$ in the singular case.

The existence of second-order derivatives requires reinforcing
regularity assumptions; at the same time it suggests reparametrizing
skewness in terms of $\delta^{(1)}=\operatorname{sign}(\delta)\delta^2$
instead of $\delta$.

\subsection{Towards a singularity-free reparametrization: Second-order
scores}\label{rep1}

Letting $\delta^{(1)}=\operatorname{sign}(\delta)\delta^2$, consider the
reparametrization $\bolds{\vartheta}^{(1)}:= (\mu^{(1)},\sigma
^{(1)}, \delta^{(1)} )\pr$. The reinforced regularity assumptions we
need (at $ \varthetab_0^{(1)}=(\mu^{(1)},\sigma^{(1)},0)\pr=
\varthetab_0$) are as follows -- recall that here we only address the
singular case under which $f$ and $\Pi$ are such that $f=g_a\in{\cal
E}_\Psi$ for some $a\in\mathcal{A}$ (see (\ref{expfam})), so that
$f$ is entirely determined by $\Pi$ and the constant $a$, and we only
need strengthening Assumption \ref{assumA2}.

\renewcommand{\theassumption}{(A2$^+$)}
\begin{assumption}\label{assumA2pl}
Same as Assumption \textup{\ref{assumA2}} but moreover \textup{(i)} the
mapping $(z, \delta)\mapsto\Pi(z,\delta)$ is twice continuously
differentiable at $(z,0)$, $z\in\R$; \textup{(ii)} denoting by $z\mapsto\dot
\psi(z) = \partial_\delta\partial_z\Pi(z,\delta)|_{\delta=0}$
the derivative of $\psi$, the quantities
%
\[
\int_{-\infty}^{\infty}\psi^2(z)z
^2 f(z)\mrmd  z \quad\mbox{and}\quad \int_{-\infty}^\infty
\bigl(2a^{-1}\dot{\psi}(z)-2\psi^2(z)\bigr)^2f(z)\mrmd
z
\]
are finite.
\end{assumption}

Assumption \ref{assumA2pl}(i) ensures the existence of the second derivative
$\partial_{\delta}^2 f_{\mu^{(1)},\sigma^{(1)},\delta}^\Pi
(x)\vert_{\varthetab_0}$ (hence, via l'Hospital's rule, that of the
simple derivative $\partial_{\delta^{(1)}} f_{{\bolds\vartheta
}^{(1)}}^\Pi(x)\vert_{\varthetab_0}$), while Assumption \ref{assumA2pl}(ii)
guarantees the finiteness of the corresponding 
covariance matrix. Assumption \ref{assumA2pl}(i) also entails $\partial
_\delta\partial_z\Pi(z,\delta)|_{\delta=0}=\partial_z\partial
_\delta\Pi(z,\delta)|_{\delta=0}$ for all $z\in\R$, so that this
mixed derivative indeed coincides with $\dot{\psi}(z)$ (see
Assumption \ref{assumA2pl}(ii)). As already pointed out, Assumption \ref
{assumA2pl} not only
reinforces Assumption \ref{assumA2} but also, via the requirement that
$f=g_a\in{\cal
E}_\Psi$ for some $a\in\mathcal{A}$, entails Assumption \ref{assumA1},
which is no
longer needed.

Under Assumption \ref{assumA2pl}, differentiating $\log f_{\mu
^{(1)},\sigma
^{(1)}, \delta^{(1)}}^\Pi$ with respect to $\delta^{(1)}$
and, at $\delta^{(1)}=\delta= 0$, applying l'Hospital's rule once
leads to (with $ f_{\mu^{(1)},\sigma^{(1)}, \delta}^\Pi(x)$ as in
(\ref{frepar1}))
%
%
\begin{equation}
\label{special} \partial_{\delta^{(1)}} \log f_{\varthetab^{(1)}}^\Pi(x)
= \cases{\displaystyle \frac{1}{2\sqrt{\vert\delta^{(1)}\vert}}\partial_\delta\log
f_{\mu^{(1)},\sigma^{(1)},
\delta}^\Pi(x)\bigg|_{\delta=\operatorname{sign}(\delta^{(1)})(\delta
^{(1)})^{1/2}}, &\quad if $\delta^{(1)}
\neq0$,
\vspace*{2pt}\cr
\displaystyle \pm\frac{1}{2}\partial_{\delta}^2 \log
f_{\mu^{(1)},\sigma^{(1)}, \delta}^\Pi(x)|_{\delta=0}, &\quad if $\delta^{(1)}
=0$,}
\end{equation}
where the undetermined sign at $\delta= 0$ is due to the fact that the
left derivative (minus sign) and the right derivative (plus sign) do
not coincide. It follows that, with the same sign indeterminacy,
%
%
\begin{equation}
\label{score1} \partial_{\delta^{(1)}} \log f_{\varthetab^{(1)}}^\Pi(x)
\vert_{\varthetab_0^{(1)}} = \pm2 \bigl[ {a^{-1}}\dot{\psi}\bigl(\sigma
^{-1}(x-\mu)\bigr) - \psi^2\bigl(\sigma^{-1} (x-
\mu)\bigr) \bigr],
\end{equation}
hence, in line with Section~\ref{mismatch},
\begin{eqnarray*}
{\bolds\ell}_{f;\varthetab^{(1)}_0}(x) &:=&
\bigl(
\ell^1_{f;\varthetab^{(1)}_0}(x), \ell^2_{f;\varthetab^{(1)}_0}(x),
\ell^3_{f;\varthetab^{(1)}_0}(x) \bigr)\pr
\\
&:=& \pmatrix{ \partial_{\mu^{(1)}} \log f_{\varthetab^{(1)}}^\Pi(x)
\vert_{\varthetab
_0}
\cr
\partial_{\sigma^{(1)}} \log f_{\varthetab^{(1)}}^\Pi(x)
\vert_{\varthetab
_0}
\cr
\partial_{\delta^{(1)}} \log
f_{\varthetab^{(1)}}^\Pi(x) \vert_{\varthetab_0}} \\
&=& \pmatrix{
\sigma^{-1}a\psi\bigl(\sigma^{-1} (x-\mu)\bigr)
\cr
\sigma^{-1} \bigl( \sigma^{-1}(x-\mu)a\psi\bigl(
\sigma^{-1}(x-\mu)\bigr)-1 \bigr)
\cr
\pm{2} \bigl[{a^{-1}}
\dot{\psi}\bigl(\sigma^{-1}(x-\mu)\bigr)-\psi^2\bigl(\sigma
^{-1} (x-\mu)\bigr) \bigr] }
\nonumber
\end{eqnarray*}
with covariance
%
%
\begin{eqnarray}
\label{info2}
{\bolds\Gamma}_{ f;\varthetab
_0^{(1)}}&:=&\sigma^{-1} \int
_{-\infty}^{\infty} {\bolds\ell}_{f;\varthetab_0^{(1)}}(x){\bolds
\ell}^{\prime}_{f;\varthetab
_0^{(1)}}(x) f\bigl(\sigma^{-1} (x-\mu)
\bigr)\mrmd x \nonumber\\[-8pt]\\[-8pt]
&=:& \pmatrix{ \gamma_{f;\varthetab_0^{(1)}}^{11}&0&0
\vspace*{1pt}\cr
0&
\gamma_{f;\varthetab_0^{(1)}}^{22}&\pm\gamma_{f;\varthetab
_0^{(1)}}^{23}
\vspace*{1pt}\cr
0&\pm\gamma_{f;\varthetab_0^{(1)}}^{23}&\gamma_{f;\varthetab
_0^{(1)}}^{33}},\nonumber
\end{eqnarray}
where (finiteness of the integrals below follows from Assumption \ref{assumA2pl}(ii))
\begin{eqnarray*}
\gamma_{f;\varthetab_0^{(1)}}^{11}&=&a^2\sigma^{-2}\int
_{-\infty
}^\infty\psi^2(z)f(z)\mrmd z,\\ 
\gamma_{f;\varthetab_0^{(1)}}^{22}&=&\sigma^{-2}\int
_{-\infty
}^\infty\bigl(a\psi(z)z-1\bigr)^2f(z)\mrmd z,
\\
\gamma_{f;\varthetab_0^{(1)}}^{33}&=&4\int_{-\infty}^\infty
\bigl(a^{-1}\dot{\psi}(z)-\psi^2(z)\bigr)^2f(z)\mrmd z
\end{eqnarray*}
and
\[
\gamma_{f;\varthetab_0^{(1)}}^{23}=2\sigma^{-1} \int
_{-\infty
}^\infty\bigl(a\psi(z)z-1\bigr)
\bigl(a^{-1}\dot{\psi}(z)-\psi^2(z)\bigr)f(z)\mrmd  z.
\]
The existence of a left and a right score for $\delta^{(1)}$ at
$\delta^{(1)}=0$, a fact which does not occur with $\frac
{1}{2}\partial_\delta^2 \log f_{\mu^{(1)},\sigma^{(1)},\delta}^\Pi
(x)\vert_{\varthetab_0}$, is not a problem, as the linear term in the
Taylor expansion of $\log f_{{\bolds\vartheta^{(1)}_0}+(0,0,\tau
_3)\pr}^\Pi(x)$ about $\log f_{{\bolds\vartheta^{(1)}_0}}^\Pi(x)$
now is of the form $
\frac{|\tau_3|}{2}\partial_\delta^2 \log f_{\mu^{(1)},\sigma
^{(1)}, \delta}^\Pi(x)\vert_{\varthetab_0}$, so that only the sign
of $\tau_3$ gets lost, as already mentioned.

Now, let us first assume that ${\bolds\Gamma}_{ f;\varthetab
_0^{(1)}}$ has full rank. Under Assumption \ref{assumA2pl} and the new
$\bolds
{\vartheta}^{(1)}$-parametrization, the model enjoys all the
properties required for traditional root-$n$ maximum likelihood and
Lagrange Multiplier or Rao score tests.
For instance, denoting by $X_1,\ldots, X_n$ an i.i.d. sample of size
$n$ from $f^\Pi_{\varthetab_0^{(1)}}$, the Lagrange Multiplier test
rejects the null hypothesis of symmetry (in favor of an asymmetry of
unspecified sign) whenever the quadratic statistic
\[
\frac{n^{-1}
( \sum_{i=1}^n (\ell^3_{f;{\hat{\varthetab}_0^{(1)}}}(X_i
) - (\gamma_{f;\hat{\varthetab} _0^{(1)}}^{23}/ \gamma_{f;\hat
{\varthetab} _0^{(1)}}^{22})\ell^2_{f;\hat{\varthetab}_0^{(1)}}(X_i )
) )^{2}}{
\gamma_{f;\hat{\varthetab} _0^{(1)}}^{33} - ( \gamma_{f;\hat
{\varthetab} _0^{(1)}}^{23})^2 / \gamma_{f;\hat{\varthetab} _0^{(1)}}^{22}
}
\]
($\hat{\varthetab} _0^{(1)}=(\hat{\mu}, \hat{\sigma}, 0)$ stands
for a root-$n$ consistent, under $\delta=0$, estimator of $\varthetab
_0^{(1)}=\varthetab_0$)
exceeds the chi-square quantile (one degree of freedom) of order
$(1-\alpha)$.

Summing up, provided that ${\bolds\Gamma}_{ f;\varthetab
_0^{(1)}}$ has full rank, root-$n$ consistency/contiguity rates are
achieved for $\delta^{(1)}=\operatorname{sign}(\delta)\delta^2$. This implies
the same root-$n$ rates at any $\delta\neq0$; at $\delta=0$,
however, an $n^{1/2}$ rate for $\delta^{(1)}$ means an $n^{1/4}$ rate
for $\delta= \operatorname{sign}(\delta^{(1)})\sqrt{\vert\delta
^{(1)}\vert
}$. Note, however, that, despite the fact that a Gram--Schmidt argument
was used in the construction, ${\bolds\Gamma}_{ f;\varthetab
_0^{(1)}}$ in general is not diagonal: there is no reason, indeed, for
the new score (\ref{score1}) to be orthogonal to those for $\mu
^{(1)}$ and $\sigma^{(1)}$.

We have assumed, so far, that ${\bolds\Gamma}_{ f;\varthetab _0^{(1)}}$
has full rank. In most cases, the components of the new score vector
$(\ell_{f;\varthetab_0^{(1)}}^{1},\ell_{f;\varthetab
_0^{(1)}}^{2},\ell_{f;\varthetab_0^{(1)}}^{3})\pr$ are not\vspace*{1pt}
collinear anymore, so that ${\bolds\Gamma}_{ f;\varthetab_0^{(1)}}$
indeed is nonsingular; our objective of a singularity-free
parametrization\vadjust{\goodbreak} thus is achieved, with consistency rate, in the
vicinity of symmetry, of $n^{1/4} $ for $\delta$. 
But this is not a general rule: in the case of the skew-normal family,
for instance, Chiogna \cite{C05} showed that the correct rate for
$\delta$
is only $n^{1/6}$. The explanation, as we shall see in the next
section, lies in a \textit{double singularity} phenomenon, which
occurs when $\ell_{f;\varthetab_0^{(1)}}^{2}$ and $\ell
_{f;\varthetab_0^{(1)}}^{3}$ in turn are collinear (by construction,
the location score $\ell_{f;\varthetab_0^{(1)}}^{1}$, at $\varthetab
_0^{(1)}$, is orthogonal to the other two scores).

\section{Double singularity}\label{step2}
\subsection{Double singularity: A special role for Gaussian
kernels}\label{double}
The double singularity phenomenon ($\ell_{f;\varthetab_0^{(1)}}^{2}$
and $\ell_{f;\varthetab_0^{(1)}}^{3}$ collinear) takes place if and
only if
\[
b \bigl(az\psi(z)-1\bigr)/\sigma=({2}/{a})\dot{\psi}(z) 
-2
\psi^2(z) \qquad\mbox{a.e.}
\]
(a.e. here and in the sequel means Lebesgue-a.e.) for some
constant $b\in\R$ and a couple $(f,\Pi)$ such that $f\in{\cal
E}_\Psi$ (see (\ref{expfam})).
Rewriting this equation as
%
%
\begin{equation}
\label{Riccati} \dot{\psi}(z)=-\frac{ab}{2\sigma}+\frac{a^2b}{2\sigma
}z\psi(z)+a
\psi^2(z) \qquad\mbox{a.e.}
\end{equation}
yields a classical Riccati equation, whose solutions are of the form
%
%
\begin{equation}
\label{admsol} \psi(z)= \frac{-ab}{2\sigma}z
\end{equation}
or
%
%
\begin{equation}
\label{Ric} \psi(z)= \frac{-ab}{2\sigma}z+\exp\biggl(- \frac
{a^2bz^2}{4\sigma
} \biggr)\bigg/
\biggl(c-a\int_0^z\exp\biggl(-\-
\frac{a^2by^2}{4\sigma
} \biggr)\mrmd y \biggr),\qquad b, c\in\R.\quad
\end{equation}
%
First, note that $b$ has to be negative, as otherwise $\varphi
_f(z)=a\psi(z)$ would tend to $-\infty$ irrespective of the sign of
$a$ when $z\rightarrow\infty$, implying positive values of $\dot
{f}$ in the right tail of $f$, which is of course impossible for a
density function. Furthermore, since both $z\mapsto a\int_0^z\exp
(- {a^2by^2}/{4\sigma} )\mrmd y$ and $\psi$ are odd, the
constant $c$ in (\ref{Ric}) has to be zero. By (\ref{expfam}), the
natural parameter space $\mathcal{A}$ for the exponential family $
{\cal E}_\Psi$ associated with the mapping $\psi$ of (\ref{Ric})
then consists of the set of values of $a$ for which the integral
\begin{eqnarray*}
\int_{-\infty}^\infty\exp\bigl(-a\Psi(z)\bigr)\mrmd z &=&\int
_{-\infty}^\infty\exp\biggl(\frac{a^2b}{4\sigma}z^2+
\log\biggl\llvert\int_0^z\exp\biggl(-
\frac{a^2by^2}{4\sigma} \biggr)\mrmd y\biggr\rrvert\biggr)\mrmd z
\\
&=&\int_{-\infty}^\infty\exp\biggl(\frac{a^2b}{4\sigma}z^2
\biggr) \biggl\llvert\int_0^z\exp\biggl(-
\frac{a^2b}{4\sigma}y^2 \biggr)\mrmd y\biggr\rrvert \mrmd z\vadjust{\goodbreak}
\end{eqnarray*}
is finite. After a change of variable involving the quantity $\sqrt
{a^2|b|/4\sigma}$, this requirement appears to be equivalent to
%
%
\begin{equation}
\label{notfinite} \int_{-\infty}^\infty\exp
\bigl(-z^2\bigr) \biggl\llvert\int_0^z
\exp\bigl(y^2\bigr)\mrmd y\biggr\rrvert \mrmd z<\infty.
\end{equation}
However, one easily can check that
$\lim_{z\rightarrow\infty}z\exp(-z^2) \llvert\int_0^z\exp
(y^2)\mrmd y\rrvert=1/2$, which means that $\exp(-z^2) \llvert\int_0^z\exp
(y^2)\mrmd y\rrvert$ behaves as $1/z$ for large values of $z$. It
follows that (\ref{notfinite}) is impossible. Hence, the natural
parameter space $\mathcal{A}$ is empty, meaning that no symmetric
kernel $f$ associated to the mapping $\psi$ of (\ref{Ric}) can yield
singular Fisher information. Therefore, the only admissible solution
to (\ref{Riccati}) is (\ref{admsol}).

This finding is quite remarkable: combined with the fact that $f\in
{\cal E}_\Psi$ (which is equivalent to $\varphi_f=a\psi$), it
implies that double singularity only can occur for symmetric kernels
$f$ such that $\varphi_f(z) = c_1z$ for some constant $c_1$ -- namely,
for Gaussian kernels. Those Gaussian kernels moreover should be
combined with a skewing function $\Pi$ such that $\psi(z)=c_2z$ for
some constant~$c_2$.

While Fisher singularity arises as a mismatch between the symmetric
kernel and the skewing function, and hence can occur with all possible
symmetric kernels, the double singularity phenomenon thus is specific
to the Gaussian kernel, hence to a well-determined subclass of
\textit{generalized skew-normal distributions} (in the sense of
Loperfido \cite{L04},
see the Introduction). This also implies that, under the assumptions
made, $n^{1/4}$ consistency rates for $\delta$ are achieved for all
other skew-symmetric families subject to Fisher singularity.

We formalize that result in the following theorem.
%
%
\begin{theor}\label{result}
Consider the skew-symmetric family defined in (\ref{HLSS}). Then:
\begin{longlist}[(ii)]
\item[(i)] under Assumptions \textup{\ref{assumA1}} and
\textup{\ref{assumA2}}, the couple $(f,\Pi)$ leads to a skew-symmetric
family subject to Fisher singularity at $\delta=0$ if and only if the
symmetric kernel $f$ is related to the skewing function $\Pi$ via the
fact that $f\in{\cal E}_\Psi$, see (\ref{expfam});
\item[(ii)] under Assumption \textup{\ref{assumA2pl}}, the couple
$(f,\Pi)$ leads to a skew-symmetric family subject to the double
singularity phenomenon if and only if the symmetric kernel $f$ is the
normal kernel $\phi$ and the skewing function $\Pi$ moreover satisfies
$\psi(z):= \partial _\delta\Pi(z,\delta)\vert_{\delta=0}=cz$ for some
real constant~$c$; the family then is a particular case of the
generalized skew-normal family (Loperfido~\cite{L04}).
\end{longlist}
\end{theor}
This theorem completely characterizes the double singularity problem,
hence complements the simple singularity characterization of Hallin and
Ley \cite{HL12}.

\subsection{A singularity-free reparametrization}\label{32}

Still inspired by Rotnitzky \textit{et al.} \cite{Ro00}, let us now proceed
with this second singularity the way we did with the first one,
producing a second, hopefully singularity-free, reparametrization.
Since the symmetric kernel $\phi$ is the only candidate for this double
singularity phenomenon, we limit ourselves to $f=\phi$. Moreover, we
know from the previous section that $\psi$ has to be of the form
$\psi(z)=c_2z$; hence, in view of the fact that $z=\varphi
_\phi(z)=a\psi(z)$, we have $c_2=1/a$. Applying the same Gram--Schmidt
idea as in Section~\ref{GS1}, but with the score for scale $\ell
_{\phi;\varthetab^{(1)}}^{2}$\vadjust{\goodbreak} substituted for the score for location,
we project $\ell_{\phi;\varthetab^{(1)}}^{3}$ onto the subspace
orthogonal to $\ell_{\phi;\varthetab^{(1)}}^{1}$ and
$\ell_{\phi;\varthetab^{(1)}}^{2}$ in the $L_2$ geometry of the
information matrix (\ref{info2}). The resulting residual score for
skewness then, as expected, is zero at $\varthetab_0^{(1)}$; indeed
\begin{eqnarray*}
&&\ell_{\phi;\varthetab_0^{(1)}}^{3}(x)-\ell_{\phi;\varthetab
_0^{(1)}}^{2} (x){
\operatorname{Cov}\bigl(\ell_{\phi;\varthetab_0^{(1)}}^{2},\ell_{\phi
;\varthetab_0^{(1)}}^{3}
\bigr)}/{\operatorname{Var}\bigl(\ell_{\phi;\varthetab
_0^{(1)}}^{2}\bigr)}
\\
&&\quad= \pm\frac{2}{a^2} \biggl(1- \frac{(x-\mu)^2}{\sigma^2} \biggr)
\\
&&\qquad{} -\sigma^{-1} \biggl(\frac{(x-\mu)^2}{\sigma
^2}-1 \biggr)
\frac{\pm2\sigma^{-1} \int_{-\infty}^\infty(z^2
-1)(a^{-2} -a^{-2}z^2)\phi(z)\mrmd  z}{\sigma^{-2}\int_{-\infty}^\infty
(z^2-1)^2\phi(z)\mrmd z}
\\
&&\quad= 0.
\end{eqnarray*}

Translating, as in Section~\ref{GS1}, this projection in terms of
parameters leads to a reparametrization of the form $(\mu^{(2)},
\sigma^{(2)},\delta)\pr$, with
\[
\mu^{(2)}=\mu^{(1)}=\mu+2\delta\sigma/{a}
\]
and
\[
\sigma^{(2)}=\sigma^{(1)}+\delta^{(1)} \frac{
\operatorname{Cov}(\ell_{\phi;\varthetab_0^{(1)}}^{2},\ell_{\phi
;\varthetab_0^{(1)}}^{3})}{
\operatorname{Var}(\ell_{\phi;\varthetab_0^{(1)}}^{2})}=\sigma
^{(1)}\bigl(1-2\delta^2/a^2
\bigr).
\]

Note that this\vspace*{1pt} reparametrization again is global, although the
orthogonality argument above only holds at $ \delta=0$. Also note that
$\sigma^{(2)}$ is not necessarily positive; although $\mu^{(2)}$ and
$\sigma^{(2)}$ jointly characterize location and scale, they cannot be
interpreted separately as a location and a scale parameter. In line
with previous notation, we denote by $f_{\mu^{(2)},\sigma
^{(2)},\delta}^\Pi$ the resulting skew-symmetric density, keeping in
mind the fact that the symmetric kernel $f$ is $\phi$. The same
remarks as for the first reparametrization are in order: keeping
$\delta$ as the skewness parameter yields $n^{1/6}$
consistency/contiguity rates. Indeed, the first two derivatives with
respect to $\delta$ now cancel at $\delta= 0$, so that derivatives of
order three play the dominant role in local approximations of log-likelihoods.

In the particular case of the skew-normal family, similar ideas have
been exploited by Chiogna \cite{C05}, where the $n^{1/6}$ rates also are
established. The reparametrization developed there, however, does not
coincide with ours, as the orthogonalization it is based on holds at
one specific value $\varthetab_0^*=(\mu^*,\sigma^*, 0)$ of
$\varthetab_0$ only ($(\mu^*,\sigma^*)$ arbitrary but fixed). The
resulting score for skewness accordingly vanishes at $\varthetab_0^*$,
while ours vanishes at all $\varthetab_0$. 

Appearance of third derivatives, in turn, suggests using $\delta
^{(2)}=\delta^3$ as a new parameter of skewness, yielding the
reparametrization ${\bolds\vartheta}^{(2)}:=(\mu^{(2)},\sigma
^{(2)},\delta^{(2)})\pr$, with
${\bolds\vartheta}^{(2)}_0:=(\mu,\sigma,0)\pr=\varthetab_0$. The new
score for skewness then will be calculated according to
%
%
\begin{equation}
\label{Hosp} \partial_{\delta^{(2)}} \log f_{\varthetab^{(2)}}^\Pi(x)
= \cases{\displaystyle \frac{1}{3(\delta^{(2)})^{2/3}}\partial_\delta\log
f_{\mu
^{(2)},\sigma^{(2)}, \delta}^\Pi(x)\bigg|_{\delta=(\delta^{(2)})^{1/3}},
&\quad
if $\delta^{(2)}
\neq0$,
\vspace*{2pt}\cr
\displaystyle \frac{1}{6}\partial_{\delta}^3 \log
f_{\mu^{(2)},\sigma^{(2)},
\delta}^\Pi(x)|_{\delta=0}, &\quad if $\delta^{(2)}
=0$,}
\end{equation}
and follows on applying l'Hospital's rule twice. This, however,
requires the following reinforcement of Assumption \ref{assumA2pl}.

\renewcommand{\theassumption}{(A2$^{++}$)}
\begin{assumption}\label{assumA2plpl}
Same as Assumption \textup{\ref{assumA2pl}}, but now \textup{(i)} the
mapping $(z,\delta)\mapsto\Pi(z,\delta)$ is three times continuously
differentiable at $(z,0)$ for all $z\in\mathbb{R}$; \textup{(ii)}
letting $\Upsilon(z):=\partial^3_\delta\Pi(z,\delta)\vert_{\delta =0}$,
the integral
$\int_{-\infty}^\infty(\frac{8}{3a^3}z^3-\frac
{8}{a^3}z+\frac{1}{3}\Upsilon(z) )^2\phi(z)\mrmd  z$ is finite.
\end{assumption}

Assumption \ref{assumA2plpl}(i) ensures\vspace*{-1pt} the existence of the third-order
derivative $\partial_\delta^3f_{\mu^{(2)},\sigma^{(2)},\delta}^\Pi
$ at $\varthetab_0^{(2)}=\varthetab_0$, while
Assumption \ref{assumA2plpl}(ii) guarantees finiteness of the corresponding
covariance matrix. Also note that the mixed derivative $\partial
_z\partial^2_\delta\Pi(z,\delta)|_{\delta=0}$ vanishes by the
definition of skewing functions; so does $\partial^2_z\partial_\delta
\Pi(z,\delta)|_{\delta=0}=\partial^2_z\psi(z)$ for all $z$, since
we are dealing (Theorem \ref{result}(ii)) with skewing functions such
that $\psi
(z)=z/a$ is linear. Finally note that $\Upsilon(z)$, by (\ref{c1}),
is an odd function. These facts greatly simplify calculations.

Assumption \ref{assumA2plpl} implies, for this second
reparametrization, the existence, at $\varthetab_0$, of a third-order
score vector ${\bolds\ell}_{\phi;\varthetab_0^{(2)}}$ with finite
covariance matrix ${\bolds\Gamma}_{\phi;\varthetab_0^{(2)} }$, enjoying
the same properties as the second-order scores described in
Section~\ref{rep1}. Elementary algebra yields
\begin{eqnarray*}
{\bolds\ell}_{\phi;\varthetab^{(2)}_0}(x)&:=& \pmatrix{ {\ell}_{\phi
;\varthetab^{(2)}_0}^1
\vspace*{1pt}\cr
{\ell}_{\phi;\varthetab^{(2)}_0}^2
\vspace*{1pt}\cr
{\ell}_{\phi;\varthetab^{(2)}_0}^3}
:=
\pmatrix{ \partial_{\mu^{(2)}} \log
f_{\varthetab^{(2)}}^\Pi(x) \vert_{\varthetab
_0}
\vspace*{1pt}\cr
\partial_{\sigma^{(2)}} \log f_{\varthetab^{(2)}}^\Pi(x) \vert
_{\varthetab
_0}
\vspace*{1pt}\cr
\partial_{\delta^{(2)}} \log
f_{\varthetab^{(2)}}^\Pi(x) \vert_{\varthetab
_0}}
\\
& = & \pmatrix{ \sigma^{-1} \bigl(\sigma^{-1}(x-\mu) \bigr)
\vspace*{1pt}\cr
\sigma^{-1} \bigl(\bigl(\sigma^{-1}(x-\mu)
\bigr)^2-1 \bigr)
\vspace*{1pt}\cr
\displaystyle \frac{8}{3a^3} \bigl(\sigma^{-1}(x-
\mu) \bigr)^3-\frac{8}{a^3}\sigma^{-1}(x-\mu)+
\frac{1}{3}\Upsilon\bigl(\sigma^{-1}(x-\mu) \bigr)}
\end{eqnarray*}
and
\[
{\bolds\Gamma}_{\phi;\varthetab_0^{(2)}}:=\sigma^{-1} \int
_{-\infty}^{\infty} {\bolds\ell}_{\phi;\varthetab
_0^{(2)}}(x){\bolds
\ell}^{\prime}_{\phi;\varthetab_0^{(2)}}(x) \phi\bigl(\sigma^{-1} (x-\mu)
\bigr)\mrmd x =: \pmatrix{ \gamma_{\phi;\varthetab_0^{(2)}}^{11}&0&
\gamma_{\phi;\varthetab
_0^{(2)}}^{13}
\cr
0&\gamma_{\phi;\varthetab_0^{(2)}}^{22}&0
\cr
\gamma_{\phi;\varthetab_0^{(2)}}^{13}&0&\gamma_{\phi;\varthetab
_0^{(2)}}^{33}}
\]
with
\begin{eqnarray*}
\gamma_{\phi;\varthetab_0^{(2)}}^{11}&=&\sigma^{-2}\int
_{-\infty
}^\infty z^2\phi(z)\mrmd z=
\sigma^{-2},\qquad 
\gamma_{\phi;\varthetab_0^{(2)}}^{22}=
\sigma^{-2}\int_{-\infty
}^\infty
\bigl(z^2-1\bigr)^2\phi(z)\mrmd z=2\sigma^{-2},
\\
\gamma_{\phi;\varthetab_0^{(2)}}^{33} &=& \int_{-\infty}^\infty
\biggl(\frac{8}{3a^3}z^3 -\frac{8}{a^3}z+\frac{1}{3}
\Upsilon(z) \biggr)^2 \phi(z)\mrmd  z
\end{eqnarray*}
and
\[
\gamma_{\phi;\varthetab_0^{(2)}}^{13} =\frac{1}{3}\sigma^{-1}
\int_{-\infty}^\infty z\Upsilon(z) \phi(z)\mrmd z.
\]
Here again, $\gamma_{\phi;\varthetab_0^{(2)}}^{13}$, in general, is
not zero, and ${\bolds\Gamma}_{\phi;\varthetab_0^{(2)}}$ is
not diagonal.

If we assume, as in Section~\ref{rep1}, that
${\bolds\Gamma}_{\phi;\varthetab_0^{(2)} }$ has full rank, denoting by
$X_1,\ldots, X_n$ an i.i.d. sample of size $n$ from $f^\Pi_{\varthetab
_0^{(2)}}$, the score vector ${\bolds\ell}_{\phi;\varthetab^{(2)}_0}$
provides a linear term to the Taylor expansion of the log-likelihood,
as well as a Lagrange Multiplier-type test of the null hypothesis of
symmetry (in the generalized skew-normal family
under study), based on the quadratic test statistic 
\[
\frac{1}{n} \sum_{i=1}^n \bigl(
\ell^3_{\phi;{\hat{\varthetab}_0^{(2)}}}(X_i ) - \sigma^{2}
\gamma_{\phi;\hat{\varthetab} _0^{(2)}}^{13}
\ell^1_{\phi;\hat{\varthetab}_0^{(2)}}(X_i
) \bigr)^{2}/ \bigl({ \gamma_{\phi;\hat{\varthetab} _0^{(2)}}^{33} -
\sigma^{2}\bigl( \gamma_{\phi;\hat{\varthetab} _0^{(2)}}^{13}
\bigr)^2 } \bigr),
\]
where $\hat{\varthetab} _0^{(2)}$ is, under the null hypothesis of
symmetry, a root-$n$ consistent estimator of location and scale. The
consistency/contiguity rate for $\delta$ (still, at $\delta=0$) is
$n^{1/6}$ while that for $\delta^{(2)}$ uniformly remains $n^{1/2}$,
and the same comments can be made as in Section~\ref{rep1}. The
particular case of the skew-normal family, from the point of view of Le
Cam's asymptotic theory of statistical experiments, is studied in full
detail by Hallin, Ley and Monti \cite{HLM12}, thereby generalizing and
extending previous work by Salvan \cite{S86} where, despite the
singularities, a locally optimal test for normality against skew-normal
alternatives is derived.

\section{Higher-order singularities}\label{step3}

It may happen, however, that ${\bolds\Gamma}_{\phi;\varthetab
_0^{(2)} }$ in turn is singular, the new third-order score for
skewness $\ell^3_{\phi;{{\varthetab}_0^{(2)}}}$ being\vspace*{1pt} (at
${\varthetab}_0^{(2)}$) collinear to the score for
location $\ell^1_{\phi;{{\varthetab}_0^{(2)}}}$ (note that, very clearly, it is orthogonal to the score for scale $\ell
^2_{\phi;{{\varthetab}_0^{(2)}}}$). If this occurs, one has to go yet
one step further with the approximation of log-likelihoods, assuming
the existence of fourth-order derivatives and ending up with $n^{1/8}$
consistency/contiguity rates for $\delta$ (but keeping $n^{1/2}$ rates
for a reparametrization $\delta^{(3)}=\operatorname{sign}(\delta)\delta
^4$ of
skewness). That $n^{1/8}$ rate, however, as we shall see, is the worst
possible one.

We are skipping details, as they are very much the same as in
Section~\ref{step2}. In order for $\ell^3_{\phi;{\varthetab
_0^{(2)}}}=\frac
{8}{3a^3}z^3-\frac{8}{a^3}z+\frac{1}{3}\Upsilon(z)$ to be linearly
dependent of $\ell^1_{\phi;{\varthetab_0^{(2)}}}=z/\sigma$,
$\Upsilon(z)$ necessarily has to be of the form $\alpha_1 z+\alpha_2
z^3$, with $\alpha_1\in\R$ and $\alpha_2=-{8}/{a^3}$ in order to
cancel the term in $z^3$. This condition on the third derivative with
respect to $\delta$ thus characterizes what we would call a \emph
{triple singularity} case (the result is formally stated in
Theorem \ref{result2} at the end of this section). It is quite easy to
construct examples suffering from this weird peculiarity; see
Section~\ref{exoho}.

The by now familiar machinery \textit{new singularity -- Gram--Schmidt
orthogonalization of scores -- reparametrization -- new higher-order
score for $\delta$} then applies,\vspace*{1pt} leading, after some direct
manipulations, to the reparametrization ${\bolds\vartheta
}^{(3)}:=(\mu^{(3)},\sigma^{(3)},\delta^{(3)})\pr$, with
\begin{eqnarray*}
\mu^{(3)}&=&\mu^{(2)}+ \biggl(-\frac{8}{a^3}+
\frac{\alpha_1}{3} \biggr)\sigma\delta^3=\mu+\frac{2}{a}\sigma
\delta+ \biggl(-\frac
{8}{a^3}+\frac{\alpha_1}{3} \biggr)\sigma
\delta^3,
\\
\sigma^{(3)}&=&\sigma^{(2)}=\sigma\bigl(1-2\delta^2/a^2
\bigr)
\end{eqnarray*}
and
\[
\delta^{(3)}=\operatorname{sign}(\delta)\delta^4.
\]
This reparametrization of skewness entails, with the same left and
right derivative interpretation for $\pm$ as in (\ref{special}),
\[
\partial_{\delta^{(3)}} \log f_{\varthetab^{(3)}}^\Pi(x) 
= \cases{\displaystyle \frac{1}{4{\vert\delta^{(3)}\vert}^{3/4}}\partial_\delta\log
f_{\mu^{(3)},\sigma^{(3)},
\delta}^\Pi(x)\bigg|_{\delta=\operatorname{sign}(\delta^{(3)})(\delta
^{(3)})^{1/4}}, &\quad if $\delta^{(3)}
\neq0$,
\vspace*{2pt}\cr
\displaystyle \pm\frac{1}{24}\partial_{\delta}^4 \log
f_{\mu^{(3)},\sigma
^{(3)}, \delta}^\Pi(x)\bigg|_{\delta=0}, &\quad if $\delta^{(3)}
=0$,}
\]
by means of a triple use of l'Hospital's rule. This, however, requires
fourth-order derivatives, hence further strengthening of
Assumption \ref{assumA2plpl}.

\renewcommand{\theassumption}{(A2$^{+++}$)}
\begin{assumption}\label{assumA2plplpl}
Same as Assumption \textup{\ref{assumA2plpl}}, but now
the mapping $(z,\delta)\mapsto\Pi(z,\delta)$ is four times
continuously differentiable at $(z,0)$, $z\in\mathbb{R}$.
\end{assumption}

Let us remark that
we do not need to assume finiteness of Fisher information for skewness,
as this, as we shall see, will always be the case after this third
reparametrization. Clearly, as in all previous cases, both the score
for location ${\ell}_{\phi;\varthetab^{(3)}_0}^1$ and the score for
scale ${\ell}_{\phi;\varthetab^{(3)}_0}^2$ remain the same as in the
original parametrization, and the new fourth-order score for skewness,
for skewing functions satisfying $\partial^3_\delta\Pi(z,\delta
)\vert_{\delta=0}=\alpha_1z-\frac{8}{a^3}z^3$, becomes (after
lengthy but quite elementary calculation)
\begin{eqnarray*}
{\ell}_{\phi;\varthetab^{(3)}_0}^3&=&\partial_{\delta^{(3)}} \log
f_{\varthetab^{(3)}}^\Pi(x) \vert_{\varthetab^{(3)} _0}
\\
&=&-\frac{10}{a^4}+\frac{2\alpha_1}{3a}+ \biggl(\frac{6}{a^4}-
\frac
{2\alpha_1}{3a} \biggr) \biggl(\frac{x-\mu}{\sigma} \biggr)^2+
\frac
{4}{3a^4} \biggl(\frac{x-\mu}{\sigma} \biggr)^4.
\end{eqnarray*}
One again easily can check that this quantity is centred under
$\varthetab_0^{(3)}$. The interesting feature, however, is that the
term $\frac{4}{3a^4} (\frac{x-\mu}{\sigma} )^4$ by no
means can cancel out, and hence prevents ${\ell}_{\phi;\varthetab
^{(3)}_0}^3$ from being any linear combination of the location and
scale scores. Thus, the resulting (at $\varthetab^{(3)}_0$) Fisher
information matrix
\begin{eqnarray*}
{\bolds\Gamma}_{\phi;\varthetab_0^{(3)}}&:=&\sigma^{-1} \int
_{-\infty}^{\infty} {\bolds\ell}_{\phi;\varthetab
_0^{(3)}}(x){\bolds
\ell}^{\prime}_{\phi;\varthetab_0^{(3)}}(x) \phi\bigl(\sigma^{-1} (x-\mu)
\bigr)\mrmd x
\\
&=& \pmatrix{ \sigma^{-2}&0&0
\vspace*{1pt}\cr
0&2\sigma^{-2}&
\displaystyle \sigma^{-1} \biggl(\frac{28}{a^4}-\frac{4\alpha
_1}{3a} \biggr)
\vspace*{2pt}\cr
0&
\displaystyle \sigma^{-1} \biggl(\frac{28}{a^4}-\frac{4\alpha_1}{3a} \biggr)&
\displaystyle \frac{1304}{3a^8}-\frac{112\alpha_1}{3a^5}+\frac{8\alpha_1^2}{9a^2}}
\end{eqnarray*}
(the finiteness of which is obvious) \emph{cannot} be singular, which
in turn implies that $n^{1/8}$ rates of convergence for $\delta$ are
the worst possible! The structural reason behind this result lies in
the fact that, by the definition of skewing functions, $\partial
^4_\delta\Pi(z,\delta)\vert_{\delta=0}$ equals zero, hence cannot
interfere in the fourth derivative, contrary to $\partial^3_\delta\Pi
(z,\delta)\vert_{\delta=0}$ which plays the crucial role in
canceling the third-order derivative.

Those results are summarized in the following theorem, which
complements Theorem~\ref{result}.

%
\begin{theor}\label{result2}
Consider the skew-symmetric family defined in (\ref{HLSS}). Then,
\begin{longlist}[(ii)]
\item[(i)] under Assumption \textup{\ref{assumA2plpl}}, the couple $(f,\Pi)$ leads
to a skew-symmetric family subject to the triple singularity phenomenon
if and only if the symmetric kernel $f$ is the normal kernel $\phi$
and the skewing function $\Pi$ moreover is such that $\psi(z):=
\partial_\delta\Pi(z,\delta)\vert_{\delta=0}=z/a$ for some
nonzero real constant $a$, and has third-order derivative $\Upsilon
(z):=\partial^3_\delta\Pi(z,\delta)\vert_{\delta=0}=\alpha
_1z-\frac{8}{a^3}z^3$ for some real constant $\alpha_1$ possibly zero.
\item[(ii)] under Assumption ($\mathrm{A}2^{+++}$), no couple $(f,\Pi)$ leads
to a skew-symmetric family subject to a fourfold/quadruple singularity
phenomenon.
\end{longlist}
\end{theor}

\section{Examples}\label{exos}

In this section, we illustrate our findings on the basis of some
well-known examples of the literature. Our presentation goes \emph
{crescendo}: starting, for the sake of completeness, with
singularity-free families, we consider simple, double, and finally
triple singularities.

\subsection{Singularity-free families}\label{exono} Famous
singularity-free examples comprise, inter alia, the
skew-exponential power distributions of Azzalini \cite{A86} with p.d.f. $2
c^{-1}\exp(-|z|^\alpha/\alpha)\Phi(\delta\operatorname
{sign}(z)|z|^{\alpha
/2}(2/\alpha)^{1/2})$ for $\alpha>1$ and $c=2\alpha^{1/\alpha
-1}\Gamma(1/\alpha)$, and the skew-$t$ distributions of Azzalini and
Capitanio \cite{AC03} with p.d.f. $2t_\nu(z)T_{\nu+1}(\delta z(\nu
+1)^{1/2}(z^2+\nu)^{-1/2})$ where $t_\eta$ and $T_\eta$,
respectively, stand for the p.d.f. and c.d.f. of a standard Student
distribution with $\eta$ degrees of freedom. These examples
are discussed at length in Hallin and Ley \cite{HL12}, where we refer
to for
details. In that same paper, an example of skewing function for which
no mismatching symmetric kernel exists is given, namely $\Pi(z,\delta
)=\Pi(\delta\sin(z))$ with $\Pi\dvtx\R\rightarrow[0,1]$ a
differentiable function satisfying $\Pi(-y)+\Pi(y)=1$ for all $y\in
\R$ and such that $\dot{\Pi}(0)=\mrmdd \Pi(y)/\mrmdd y\vert_{y=0}$ exists and
differs from zero.

\subsection{Simple singularities} As shown in Hallin and Ley \cite{HL12},
the easiest-to-construct mismatching skewing function for a given
symmetric kernel $f$ is of the form $\Pi(\delta\varphi_f(z))$, with
$\Pi$ as described above. For any symmetric kernel $f$, it is readily
seen that the location and skewness scores then are collinear.

Under the assumptions made, double singularity requires the additional
assumption that $\ddot{\Pi}(0):=d^2\Pi(y)/(\mrmdd y)^2\vert_{y=0}$ exists
and, by construction, equals zero. Theorem \ref{result} then tells us
that among the p.d.f.s $2f(z)\Pi(\delta\varphi_f(z))$ only the
skew-normal, obtained for $f=\phi$, suffers from the double
singularity. Thus all non-Gaussian kernels $f$ yield examples of simple
singularities.

\subsection{Double singularities}Concerning double singularity, a
prominent example is of course Azzalini's skew-normal family, with
p.d.f.
$2\phi(z)\Phi(\delta z)$. Let us briefly show that higher-order
singularities are excluded in that family. Straightforward calculation
yields $a=\sqrt{2\pi}$ and $\Upsilon(z)=-(2\pi)^{-1/2}z^3$, which
is different from $-\frac{8}{a^3}z^3=-(2/\pi)^{3/2}z^3$. Hence,
Theorem \ref{result2} readily yields the well-known result of
$n^{1/6}$ consistency rates for $\delta$ in the skew-normal
distribution. For the sake of completeness, we also provide for this
famous example the corresponding score for skewness, which equals\vspace*{-2.5pt}
\[
\frac{4-\pi}{3\pi\sqrt{2\pi}}z^3-\frac{4}{\pi\sqrt{2\pi}}z.
\]

Nadarajah and Kotz \cite{NK03} propose another family of skew densities
generated by the normal kernel, with p.d.f.s of the form $2\phi
(z)G(\delta z)$ where $G$ is some univariate symmetric c.d.f. They call
\emph{skew normal-$G$} the resulting families of densities. Their
definition includes as particular cases the skew normal--normal model,
the skew normal-$t$, the skew normal-Cauchy, the skew normal-Laplace,
the skew normal-logistic and the skew normal-uniform families.
Theorem \ref{result} tells us that all skew normal-$G$ models suffer
from double singularity, a fact that, except of course for the skew
normal--normal (which, up to an additional scale parameter, coincides
with the classical skew-normal), has never been noticed. Consequently,
these models have to be treated with much care when used for
inferential purposes. The problem with those families obviously stems
from the product $\delta z$ inside $G$; see Section~\ref{fc} for
further discussion of such skewing functions.

\subsection{Higher-order singularities}\label{exoho}

Let us further analyze the families of Nadarajah and Kotz \cite{NK03}.
Assume that $G$ is three times continuously differentiable. Elementary
calculations show that $a=1/g(0)$, where $g(z):=\mrmdd G(z)/\mrmdd z$, and
$\Upsilon(z)=\ddot{g}(0)z^3$. We know from Theorem \ref{result2}
that a triple singularity can only occur if $\ddot{g}(0)=-\frac
{8}{a^3}=-8(g(0))^3$. Among the distributions considered by Nadarajah
and Kotz \cite{NK03}, this equality holds for the skew normal-logistic only,
for which $g(0)=1/4$ and $\ddot{g}(0)=-1/8$. Thus, while all their
other skew normal-$G$ distributions yield $n^{1/6}$ consistency rates
for $\delta$, the skew normal-logistic one requires the worst possible
rates, namely $n^{1/8}$ rates.

Finally, consider the ``lifted'' skew-normal distribution, with p.d.f.\vspace*{-2.5pt}
%
%
\begin{equation}
\label{lift} 2\phi(z)\Phi\bigl(\delta z-(4-\pi) (6\pi)^{-1}
\delta^3z^3\bigr).
\end{equation}
Here, $a=\sqrt{2\pi}$ and $\Upsilon(z)=-(2/\pi)^{3/2}z^3=-\frac
{8}{(\sqrt{2\pi})^3}z^3=-\frac{8}{a^3}z^3$, entailing, by
Theorem \ref{result2}, a~triple singularity and hence $n^{1/8}$
consistency rates for $\delta$.
Note that this distribution is part of the so-called class of \emph
{flexible generalized skew-normal distributions} defined in Ma and
Genton \cite{MG04}. More generally, in that paper, the authors have proposed
\emph{flexible skew-symmetric distributions} with skewing functions of
the form
$\Pi(z,\delta):=\Pi(H_\ell(\delta z))$, with $\Pi$ as defined in
Section~\ref{exono} and $H_\ell$ an odd polynomial of order $\ell$
(meaning that the polynomial only contains odd power terms). Since, in
the first four derivatives, all terms of the form $(\delta z)^s$ with
odd $s\geq5$ do not play any role, one can directly construct an
infinity of flexible generalized skew-normal distributions suffering
from triple singularity: take any odd polynomial $H_\ell$ with the
terms in $\delta z$ and $(\delta z)^3$ as\vadjust{\goodbreak} in (\ref{lift}), for
instance,
\[
2\phi(z)\Phi\Biggl(\delta z-(4-\pi) (6\pi)^{-1}\delta^3z^3+
\sum_{i=2}^\ell\alpha_{2i+1}(\delta
z)^{2i+1}\Biggr)
\]
with $\alpha_i\in\R$ and $2\leq\ell\in\N$.

\section{Some concluding remarks}\label{final}

We conclude this paper by a short discussion of two structural issues:
the centred parametrization (Section~\ref{cp}) and the type of
skewing function that causes most of the trouble when using a Gaussian
kernel (Section~\ref{fc}).

\subsection{The centred parametrization}\label{cp}

In order to remedy Fisher singularity problems, Azzalini \cite{A85}, in the
very same paper where he first introduces the skew-normal densities,
and for the specific case of the skew-normal family, suggested an
alternative parametrization, the so-called \textit{centred
parametrization}. Denoting by $Z$ a random variable with skew-normal
density $2\phi(z)\Phi(\delta z)$, let $Y:=\mu+\sigma Z$ a.s.
($\sigma>0$): $Y$ then has skew-normal density $2\sigma^{-1}\phi
(\sigma^{-1}(z-\mu))\Phi(\delta\sigma^{-1}(z-\mu))$. That density
has finite third-order moments: letting $\theta_1:=\mathrm{E}[Y]$ and
$\theta_2:=\mathrm{Var}^{1/2}[Y]$,
define $\gamma_1:=\mathrm{E}[(Y-\theta_1)^3]/\theta_2^3$ as $Y$'s
(hence also $Z$'s) third standardized cumulant. The triple
$\thetab:=(\theta_1,\theta_2,\gamma_1)\pr$ provides a parametrization
of the skew-normal family, the \textit{centred parametrization}
(hereafter CP). The terminology ``centred'' refers to the fact that
the new location and scale parameters $\theta_1$ and $\theta_2$ are
such that $(Y-\theta_1)/\theta_2 = (Z-\mathrm{E}[Z])/\operatorname
{Var}^{1/2}[Z]$
has mean zero and variance one, whereas the original location and scale
$\mu$ and $\sigma$ values lead to $(Y-\mu)/\sigma=Z$, which is not
centred about its mean. Azzalini calls $\mu$, $\sigma$ and $\delta$
\emph{direct parameters}, since they can be directly read from the
density of $Y$.

Besides Azzalini \cite{A85}, the centred parametrization has been
discussed in Azzalini and Capitanio~\cite{AC99}, Pewsey \cite{P00} and
Chiogna \cite{C05}, to cite but these, for the skew-normal family,
extended to the multinormal setup in Arellano-Valle and Azzalini
\cite{AA08}, and to the skew-$t$ distributions in DiCiccio and Monti
\cite{DM11} and Arellano-Valle and Azzalini \cite{AA11}.

The CP does not suffer from the Fisher singularity problem, and
provides parameters that can all be estimated at the usual $n^{1/2}$
rate. It is well-suited for inferential purposes (see, e.g.,
Pewsey~\cite{P00}), and enjoys a simple traditional moment-based
interpretation (which is the main motivation for Arellano-Valle and
Azzalini \cite{AA11} to extend it to the skew-$t$ context although
skew-$t$ families do not exhibit any Fisher singularity).

The main drawback of the CP lies in its complicated analytical form.
Expressing the centred parameters $\thetab$ in terms of the original
ones $\varthetab$ yields (for the skew-normal family)
\[
\theta_1=\mu+\sigma\sqrt{2/\pi} {\delta}\bigl(1+\delta^2
\bigr)^{-1/2},\qquad \theta_2=\sigma\bigl(1+\delta^2(1-2/
\pi)\bigr)^{1/2}\bigl(1+\delta^2\bigr)^{-1/2}
\]
and
\[
\gamma_1=\frac{4-\pi}{2} \biggl(\frac{2}{\pi}
\biggr)^{3/2}\delta^3\bigl(1+\delta^2(1-2/\pi)
\bigr)^{-3/2}.
\]
The success of skew-symmetric families is largely due to the analytical
simplicity of the decomposition (\ref{HLSS}) of a skew density into
the product of a symmetric kernel and a skewing function. That
simplicity is closely related to the original $\varthetab$
parametrization, and gets lost in the CP. So is the flexibility that
allows for combining various symmetric kernels and skewing functions
in (\ref{HLSS})\vadjust{\goodbreak} without altering its structure: contrary to the
mapping $\varthetab\mapsto f_{\varthetab}^\Pi(x)$, which does not
depend on $(f,\Pi)$, the mapping $\thetab\mapsto f_{\thetab}^\Pi
(x)$ very much does. Although simple from the point of view of
interpretation, the CP thus does not avoid analytical complexity, which
results into lengthy and very tedious calculations -- see
Chiogna \cite{C05}.\footnote{We take this opportunity to indicate two
unfortunate typos in that paper, namely (i) an exponent $^{-1/2}$ is
missing in the expression of $\delta$ ($\lambda$ in the notation of
the paper) in terms of $\gamma_1$ on page 338, and (ii) the score
function $u^\Theta(\phi^*,\gamma_1^*)$ on page 339, instead of three
linear terms, should involve a quadratic term for the $\theta_2$-score
and a term of the form $(x^3-3x)$ for the $\gamma_1$-score.}

Both the centred parametrization and those we are proposing in this
paper are motivated by Fisher singularity problems at $\delta=0$, and
both are losing some of the simplicity of the original parametrization.
The CP, however, is guided by interpretability considerations (the
centred parameters \textit{always} are the mean, the standard error,
and the third-order cumulant); although solving the Fisher singularity
problem in the skew-normal family, there is no guarantee it does so in
other skew-symmetric families. Our reparametrizations, on the contrary,
are guided by Fisher information considerations, and are specifically
designed to solve the Fisher singularity problem, irrespective of the
skew-symmetric family under study -- with distinctive forms for simple,
double, and triple singularity. The way they deal with $\delta$
(exponentiating it into $\operatorname{sign}(\delta)\delta^2$, $\delta
^3$, or
$\operatorname{sign}(\delta)\delta^4$) preserves its interpretation,
regardless of the chosen symmetric kernel, as a tuning quantity in the
skewing mechanism that generates the family. Finally, due to their
Gram--Schmidt nature, our reparametrizations are tailored for the
construction of optimal tests of symmetry of the Lagrange Multiplier or
Rao score type.

For the sake of comparison, we provide, in the
\hyperref[app]{Appendix}, the scores for skewness associated with each
of the two reparametrizations (the CP and ours) in the skew-normal
family.

\subsection{A brief discussion of skewing functions of the form
\texorpdfstring{$\Pi(z,\delta)=\Pi(\delta z)$}{Pi(z,delta) = Pi(delta z)}}\label{fc}

We conclude this section with a few comments on the most frequent type
of skewing function, namely $\Pi(z,\delta)=\Pi(\delta z)$ with
$\Pi\dvtx\R\rightarrow[0,1]$ satisfying $\Pi(-y)+\Pi(y)=1$ for all
$y\in\R
$ (and the required differentiability conditions). Such functions are
the most natural examples of a skewing function such that $\psi(z)$ is
linear, yielding a risky combination with the Gaussian kernel $\phi$.

The original skew-normal family of Azzalini \cite{A85} is based on $\Pi
(z,\delta) =\Phi(\delta z)$; the same type of skewing function has
been used, inter alia, by:
\begin{enumerate}[-]
\item[-] Azzalini and Capitanio \cite{AC99} for skew-symmetric
densities of
the form $2f(z)G(\delta z)$, with $G$ some univariate symmetric
distribution function (in fact, Azzalini and Capitanio proposed
multivariate skew-elliptical distributions, but elliptical symmetry
here boils down to plain univariate symmetry);
\item[-] Gupta, Chang and Huang \cite{GCH02} for their skew-uniform, skew-$t$,
skew-Cauchy, skew-Laplace and skew-logistic distributions, which all
are special cases of Azzalini and Capitanio's \cite{AC99} construction;
\item[-] Nadarajah and Kotz \cite{NK03} for their skew normal-$G$
distributions, as described in the previous sections; and by
\item[-] G{\'o}mez, Venegas and Bolfarine \cite{GTB07} for their skew $g$-normal
densities $2g(z)\Phi(\delta z)$ where, contrary to the skew normal-$G$
distributions, normality is present in the skewing function and not in
the symmetric kernel.
\end{enumerate}
As shown in this paper, skewing functions of the form $\Pi(\delta z)$
are harmless whenever the symmetric kernel is not Gaussian. In view of
this, the skew $g$-normal distributions (free of any singularity except
for $g=\phi$) are inferentially preferable to the skew normal-$G$ ones
(which exhibit at least double singularity). An early important warning
on the combination of a Gaussian kernel with such skewing functions has
been given by Pewsey \cite{P06}, who has shown that all densities of the
form $2\phi(z) G(\delta z)$ ($G$ some symmetric univariate c.d.f.) suffer
from the singularity problem. His results are thus in total agreement
with our general findings. The peculiarities of the skew-normal
distribution, which belongs to all of the above-cited classes of
distributions, have been discussed at length in the literature. We hope
that this paper sheds some more light on the structural reasons behind
those peculiarities, and provides further warning about the dangers of
Gaussian kernels in combination with skewing functions of the form $\Pi
(\delta z)$.

\begin{appendix}\label{app}

\section*{Appendix: Expressions of the score functions for the
skew-normal distribution}

In this Appendix, we provide the explicit expressions of the score
functions for skewness in the skew-normal case (at any value of the
skewness parameter, not only in the vicinity of symmetry) for both the
centred parametrization and ours (as described in Section~\ref{32}).

In our reparametrization $\varthetab^{(2)}$, the score for skewness
$\partial_{\delta^{(2)}}\log f^{\Pi}_{\varthetab^{(2)}}(x)$ takes on
the guise of a ratio $ {h_1(\mu^{(2)},\sigma^{(2)},\delta
^{(2)})}/{h_2(\mu^{(2)},\sigma^{(2)},\delta^{(2)})}$, with
\begin{eqnarray*}
&&
h_1(\mu^{(2)},\sigma^{(2)},\delta^{(2)})\\
&&\quad:=\exp
\biggl(-\frac{(\delta^{(2)})^{2/3}(\sqrt{2\pi}(\delta
^{(2)})^{1/3}\sigma^{(2)}-(\delta^{(2)})^{2/3}(x-\mu^{(2)})+\pi
(x-\mu^{(2)}))^2}{2\pi^2(\sigma^{(2)})^2} \biggr)\bigl(
\sigma^{(2)}\bigr)^{-1}
\\
&&\qquad{}\times\bigl(\pi-\bigl(\delta^{(2)}
\bigr)^{2/3}\bigr) \bigl(4\pi\bigl(\delta^{(2)}
\bigr)^{1/3}\sigma^{(2)}+\sqrt{2}\pi^{3/2}\bigl(x-
\mu^{(2)}\bigr) -3\sqrt{2\pi}\bigl(\delta^{(2)}
\bigr)^{2/3} \bigl(x-\mu^{(2)}\bigr) \bigr)
\\
&&\qquad{}-4\pi^{2}\bigl(\delta^{(2)}\bigr)^{1/3} \Phi
\biggl(\frac{\sqrt{2\pi
}(\delta^{(2)})^{2/3}\sigma^{(2)}-\delta^{(2)}(x-\mu^{(2)})+\pi
(\delta^{(2)})^{1/3}(x-\mu^{(2)})}{\pi\sigma^{(2)}} \biggr)
\\
&&\qquad{}+\frac{2}{(\sigma^{(2)})^2}\bigl(\pi-\bigl(\delta^{(2)}\bigr)^{2/3}
\bigr) \Phi\biggl(\frac{\sqrt{2\pi}(\delta^{(2)})^{2/3}\sigma
^{(2)}-\delta
^{(2)}(x-\mu^{(2)})+\pi(\delta^{(2)})^{1/3}
(x-\mu^{(2)})}{\pi\sigma^{(2)}} \biggr)
\\
&&\qquad{} \times\bigl(-\sqrt{2}\pi^{3/2}\sigma^{(2)}
\bigl(x-\mu^{(2)}\bigr)-2\delta^{(2)}\bigl(x-\mu^{(2)}
\bigr)^2+3\sqrt{2\pi}\bigl(\delta^{(2)}\bigr)^{2/3}
\sigma^{(2)}\bigl(x-\mu^{(2)}\bigr)
\\
&&\qquad{} + 2\pi\bigl(\delta^{(2)}\bigr)^{1/3}\bigl(
\bigl(x-\mu^{(2)}\bigr)^2-\bigl(\sigma^{(2)}
\bigr)^2\bigr)\bigr)
\end{eqnarray*}
and
$ {h_2(\mu^{(2)},\sigma^{(2)},\delta^{(2)}):= 6\pi
^{2}(\pi-(\delta^{(2)})^{2/3})(\delta^{(2)})^{2/3}}$
\[
\times\Phi\biggl(\frac{\sqrt{2\pi}(\delta
^{(2)})^{2/3}\sigma^{(2)}-\delta^{(2)}(x-\mu^{(2)})+\pi(\delta
^{(2)})^{1/3}(x-\mu^{(2)})}{\pi\sigma^{(2)}} \biggr).
\]

In the CP reparametrization $\thetab=(\theta_1,\theta_2,\gamma_1)\pr$,
the score for skewness $\partial_{\gamma_1}\log f^{\Pi
}_{(\theta_1,\theta_2,\gamma_1)}(x)$ takes on the guise
$h^{\mathrm{CP}}_1(\theta_1,\theta_2,\gamma_1)+
{h^{\mathrm{CP}}_2(\theta_1,\theta
_2,\gamma_1)}/{h^{\mathrm{CP}}_3(\theta_1,\theta_2,\gamma_1)}$, with
\begin{eqnarray*}
&&h^{\mathrm{CP}}_1(\theta_1,
\theta_2,\gamma_1)=-\frac{1}{3\gamma
_1^{1/3} (\gamma_1^{2/3}+ (({4-\pi})/{2} )^{2/3}
)}\\
&&\hphantom{h^{\mathrm{CP}}_1(\theta_1,
\theta_2,\gamma_1)=}
{}+
\frac{1}{3} \biggl(\frac{2}{4-\pi} \biggr)^{2/3}
\frac{
(x-\theta_1+ ({2}/({4-\pi}) )^{1/3}\gamma
_1^{1/3}\theta_2 )^2}{\gamma_1^{1/3}\theta_2^2 (1+\gamma
_1^{2/3} ({2}/({4-\pi}) )^{2/3} )^2}
\\
&&
\hphantom{h^{\mathrm{CP}}_1(\theta_1,
\theta_2,\gamma_1)=}
{} -\frac{1}{3} \biggl(\frac{2}{4-\pi}
\biggr)^{1/3}\frac{x-\theta_1+ ({2}/({4-\pi}) )^{1/3}\gamma
_1^{1/3}\theta_2}{\gamma_1^{2/3}\theta_2 (1+\gamma
_1^{2/3} ({2}/({4-\pi}) )^{2/3} )},
\\
&&h^{\mathrm{CP}}_2(\theta_1,\theta_2,\gamma_1)\\
&&\quad=2^{1/6}\exp\biggl(\frac{\pi(4-\pi)^{-2/3}\gamma_1^{2/3}
(x-\theta_1+ ({2}/({4-\pi}) )^{1/3}\gamma
_1^{1/3}\theta_2 )^2}{2^{4/3} (1+ ({2}/({4-\pi
}) )^{2/3}\gamma_1^{2/3} ) (-1+2^{-1/3}(\pi-2)(4-\pi
)^{-2/3}\gamma_1^{2/3} )\theta_2^2} \biggr)
\\
&&\qquad{} \times\bigl((x-\theta_1) \bigl(2(4-
\pi)^{2/3}(\pi-2) \gamma_1^{4/3} +2^{2/3}(
\pi-4)^2 \bigr)\\
&&\qquad\quad\hspace*{4pt}{}+\theta_2 \bigl(2^{2/3}(4-
\pi)^2 \gamma_1 +4(4-\pi)^{5/3}
\gamma_1^{1/3} \bigr) \bigr)
\end{eqnarray*}
and
\begin{eqnarray*}
&&h^{\mathrm{CP}}_3(\theta_1,\theta_2,\gamma_1)\\
&&\quad=3(4-\pi)^{7/3}\theta_2\gamma_1^{2/3}
\biggl(1+\gamma_1^{2/3} \biggl(\frac{2}{4-\pi}
\biggr)^{2/3} \biggr)^{3/2} \bigl(2-2^{2/3}(\pi-2) (4-
\pi)^{-2/3}\gamma_1^{2/3} \bigr)^{3/2}
\\
&&\qquad{} \times\Phi\biggl(\frac{2^{1/3}\sqrt{\pi}(4-\pi
)^{-1/3}\gamma_1^{1/3} (x-\theta_1+\gamma_1^{1/3}\theta_2
({2}/({4-\pi}) )^{1/3} )}{\theta_2
(1+2^{2/3}(4-\pi)^{-2/3}\gamma_1^{2/3} )^{1/2}
(2-2^{2/3}(\pi-2)(4-\pi)^{-2/3}\gamma_1^{2/3} )^{1/2}} \biggr).
\end{eqnarray*}

Both score functions look equally complex. Under symmetry (either
$\delta^{(2)}=0$ or $\gamma_1=0$), they both yield an indetermination
of the form $0/0$, apparently calling for tedious applications of
l'Hospital's rule. For the $\varthetab^{(2)}$ reparametrization we are
proposing, such algebra is not required, though, as explicit
expressions involving higher order derivatives already have been
derived analytically (and follow quite directly from the Gram--Schmidt
orthogonalization structure): see~(\ref{Hosp}). In the CP
reparametrization $\thetab$, the required algebra is so tedious as to
defeat our version of \textsc{Mathematica}.
\end{appendix}

\section*{Acknowledgements}

Marc Hallin is member of the Acad\'emie Royale de Belgique and ECORE,
and an extra-muros Fellow of CentER, Tilburg University. His research
is supported by the Sonderforschungsbereich ``Statistical modeling of
nonlinear dynamic processes'' (SFB 823) of the German Research
Foundation (Deutsche Forschungsgemeinschaft), the Belgian Science
Policy Office (2012--2017) Interuniversity Attraction Poles, and a Discovery Grant of
the Australian Research Council.

Christophe Ley thanks the Fonds National de la Recherche Scientifique,
Communaut\'{e} fran\c{c}aise de Belgique, for support via a Mandat de
Charg\'{e} de Recherche FNRS.

Both authors would like to thank two anonymous referees for their
helpful comments that led to an improvement of the present paper, as
well as Adelchi Azzalini and Reinaldo Arellano-Valle for interesting
discussions on the centred parametrization.


%

\printhistory

\end{document}